\documentclass[a4paper,11pt,reqno]{amsart}
\usepackage{amssymb}
\usepackage{amsmath}
\usepackage{amsbsy}
\usepackage{amsthm}
\usepackage{mathtools}

\usepackage{amssymb}
\usepackage{latexsym}
\usepackage{amsmath}
\usepackage{amsbsy}
\usepackage{amsthm}
\usepackage{txfonts}
\usepackage{pxfonts}
\usepackage{pict2e}
\usepackage{mathtools}
\usepackage{amscd}
\usepackage{mathrsfs}

\newtheorem{defi}{Definition}
\newtheorem{theo}{Theorem}

\newtheorem{prop}{Proposition}[section]
\newtheorem{cor}{Corollary}[prop]

\newcommand{\Lan}{\mbox{$\mathcal{L}$}}
\newcommand{\LanP}{\mbox{$\mathcal{L}^{\ast}$}}

\newcommand{\fA}{\mbox{$\mathfrak{A}$}}

\newcommand{\fB}{\mbox{$\mathfrak{B}$}}
\newcommand{\fC}{\mbox{$\mathfrak{C}$}}

\newcommand{\fG}{\mbox{$\mathfrak{G}$}}
\newcommand{\fM}{\mbox{$\mathfrak{M}$}}

\newcommand{\fW}{\mbox{$\mathfrak{W}$}}
\newcommand{\Fil}{\mbox{$\mathcal{F}$}}
\newcommand{\Mat}{\mbox{$\mathcal{M}$}}
\newcommand{\Alg}{\mbox{$\mathfrak{F}_{\mathcal{L}}$}}
\newcommand{\AlgS}{\mbox{$\mathfrak{F}_{\mathcal{S}}$}}
\newcommand{\Thy}{\mbox{$\Sigma_{\mathcal{S}}$}}
\newcommand{\Sr}{\mbox{$\mathcal{S}_{u}$}}

\newcommand{\cS}{\mbox{$\mathscr{S}$}}
\newcommand{\cmS}{\mbox{$\mathcal{S}$}}
\newcommand{\deduc}{\mbox{$\vdash_\mathcal{S}$}}
\newcommand{\modelM}{\mbox{$\models_{\mathcal{M}}$}}
\newcommand{\cn}{\mbox{\textit{Cn}}}
\newcommand{\Cn}[1]{\mbox{$\textit{Cn}_{\mathcal{S}}(#1)$}}
\newcommand{\g}{\textit{\textbf{g}}}

\newcommand{\Var}{\mathcal{V}}
\newcommand{\VarP}{\mathcal{V}_{\mathcal{L}}}
\newcommand{\VarPP}{\mathcal{V}_{\mathcal{L}^{\ast}}}
\newcommand{\Fr}{\mbox{$\textit{\textbf{Fr}}_{\mathcal{L}}$}}
\newcommand{\Th}{\mbox{$\textit{\textbf{T}}_{\mathcal{S}}$}}

\newcommand{\fsub}{\mbox{$\Subset$}}

\newcommand{\ra}{\mbox{$\rightarrow$}}
\newcommand{\con}{\mbox{$\wedge$}}
\newcommand{\dis}{\mbox{$\vee$}}

\newcommand{\Ga}{\mbox{$\Gamma$}}

\newcommand{\on}{\mbox{$\mathbf{1}$}}

\newcommand{\set}[2]{\{#1 ~|~#2\}}

\begin{document}

\renewcommand{\thefootnote}{$\star$}

\title[Lindenbaum Method]{Lindenbaum method (propositional language)$^\star$}
\footnotetext{This is a slightly revised version of an article under the same title  written in 2013 for  \textit{Encyclopedia of Mathematics}, currently existing  (possibly in an updated form) at URL:http://www.encyclopediaofmath.org/index.php?title=Lindenbaum{\_}method{\&}oldid=30629.}
\author{Alex Citkin}
\author{Alexei Muravitsky}
		
\maketitle


\noindent Lindenbaum method is named after the Polish logician Adolf Lindenbaum who prematurely and without a clear trace disappeared in the turmoil of the Second World War at the age of about 37. (Cf.~\cite{sur82}.) The method is based on the symbolic nature of formalized languages of deductive systems and opens a gate for applications of algebra to logic and, thereby, to \textit{Abstract algebraic logic}.\\

\noindent\Large{\textbf{Lindenbaum's Theorem}}\\

\noindent A formal propositional language, say {\Lan}, is understood as a nonempty set $\VarP$ of symbols $p_0, p_1,$ $\ldots p_{\gamma},\ldots$  called propositional variables and a finite set $\Pi$ of symbols $F_0, F_1,\ldots, F_n$ called logical connectives. By $\overline{\overline{\VarP}}$ we denote the cardinality of $\VarP$. For each connective $F_i$, there is a natural number $\#(F_i)$ called the arity of the connective $F_i$. The notion of a statement (or a formula) is defined as follows:
        \[
        \begin{array}{cl}
        (f_1) &\mbox{Each variable $p\in\VarP$ is a formula};\\
        (f_2) &\mbox{If $F_i$ is a connective of the arity 0, then $F_i$}\\ 
        &\mbox{is a formula};\\
        (f_3) &\mbox{If $A_1, A_2,\ldots, A_n$, $n\geq 1$, are formulas, and $F_i$}\\ 
        &\mbox{is a connective of arity $n$, then the symbolic}\\ &\mbox{expression $F_{n}A_{1}A_{2}\ldots A_n$ is a formula};\\
        (f_4) &\mbox{A formula can be constructed only according}\\ 
        &\mbox{to the rules $(f_1)-(f_3)$.}\\
        \end{array}
        \]

\noindent The set of formulas will be denoted by {\Fr} and $\mathcal{P}(\Fr)$ denotes the power set of {\Fr}. Given a set $X\subseteq\Fr$, we denote by $\Var(X)$ the set of all propositional variables that occur in the formulas of $X$. Two formulas are counted equal if they are represented by two copies of the same string of symbols. (This is the key observation on which Theorem~\ref{P:absolutely-free} is grounded.) Another key observation (due to Lindenbaum) is that {\Fr} along with the connectives $\Pi$ can be regarded as an algebra of the similarity type associated with {\Lan}, which exemplifies an \Lan-\textit{algebra}. We denote this algebra by $\Alg$. The importance of {\Alg} can already be seen from the following statement.

\begin{theo}\label{P:absolutely-free}
Algebra $\Alg$ is a free algebra of rank $\overline{\overline{\VarP}}$ with free generators $\VarP$ in the class $($variety$)$ of all {\Lan}-algebras. In other words, $\Alg$ is an absolutely free algebra of this class. {\em (Cf.~\cite{mmt87}, section 4.11)}.
\end{theo}

\noindent A useful feature of the set {\Fr} is that it is closed under (\textit{simultaneous}) \textit{substitution}. More than that, any substitution $\sigma$ is an endomorphism
\[\sigma: \Alg\longrightarrow \Alg.\]

\noindent A \textit{monotone deductive system} (or a \textit{deductive system} or simply a \textit{system}) is a relation between subsets and elements of {\Fr}. Each such system $\deduc$ is subject to the following conditions: For all $X,Y\subseteq\Fr$,
\[
\begin{array}{cl}
(s_1) &\mbox{if $A\in X$, then $X\deduc A$};\\
(s_2) &\mbox{if $X\deduc B$ for all $B\in Y$, and $Y\deduc A$, then $X\deduc A$};\\
(s_3) &\mbox{if  $X\deduc A$, then for every substitution $\sigma$, $\sigma[X]\deduc\sigma(A)$}.
\end{array}
\]
If $A$ is a formula and $\sigma$ is a substitution, $\sigma(A)$ is called a \textit{substitution instance} of $A$. Thus, by $\sigma[X]$ above, one means the set of the instances of the formulas of $X$ with respect to $\sigma$.\\

\noindent Given two sets $Y$ and $X$, we write
\[Y\fsub X\]
\noindent if $Y$ is a finite (maybe empty) subset of $X$.\\

\noindent A deductive system is said to be \textit{finitary} if, in addition, it satisfies the following:
\[
\begin{array}{cl}
~(s_4) &\mbox{if $X\deduc A$, then there is $Y\fsub X$ such that $Y\deduc A$}.
\end{array}
\]

\noindent We note that the monotonicity property
\[
\mbox{if $X\subseteq Y$ and $X\deduc A$, then $Y\deduc A$}
\]
is not postulated, because it follows from $(s_1)$ and $(s_2)$.\\

\noindent Each deductive system $\deduc$ induces a (\textit{monotone structural}) \textit{consequence operator} $\cn_{\mathcal{S}}$ defined on the power set of {\Fr} as follows: For every $X\subseteq\Fr$,

\begin{equation}\label{E:cn-system}
A\in\Cn{X} \Longleftrightarrow X\deduc A,
%\tag{$1^{\star}$}
\end{equation}
so that the following conditions are fulfilled: For all $X,Y\subseteq\Fr$ and any substitution $\sigma$,
\[
\begin{array}{cl}
(c_1) &X\subseteq\Cn{X};\quad(\textit{reflexivity})\\
(c_2) &\Cn{\Cn{X}}=\Cn{X};~(\textit{idempotency})\\
(c_3) &\mbox{if $X\subseteq Y$, then $\Cn{X}\subseteq\Cn{Y}$};~(\mbox{\textit{monotonicity}})\\
(c_4) &\sigma[\Cn{X}]\subseteq\Cn{\sigma[X]}.~(\mbox{\textit{structurality} or}\\
 &\quad\quad\quad\quad\quad\quad\quad\quad\quad\quad
 \hspace{0.15in}\mbox{\textit{substitution invariance}})
\end{array}
\]

\noindent If $\deduc$ is finitary, then
\[
\begin{array}{cl}
\hspace{-1.4in}(c_5) &\Cn{X}=\bigcup\set{\Cn{Y}}{Y\fsub X},
\end{array}
\]
in which case $\cn_{\mathcal{S}}$ is called \textit{finitary}.\\

\noindent Conversely, if an operator $\textit{Cn}:\mathcal{P}(\Fr)\rightarrow \mathcal{P}(\Fr)$ satisfies the conditions $(c_1)-(c_4)$ (with \textit{Cn} instead of $\textit{Cn}_{\mathcal{S}}$), then the equivalence
\[
X\deduc A\Longleftrightarrow A\in\textit{Cn}(X)
\]
defines a deductive system, $\mathcal{S}$. Thus
(\ref{E:cn-system}) allows one to use the deductive system and consequence operator (in a fixed formal language) interchangeably or even in one and the same context. For instance, we call 
\[\Th=\Cn{\emptyset}\] 
the set of \textit{theorems} of the system $\vdash_{\cmS}$ (i.e. {\cmS}-\textit{theorems}), and given a subset $X\subseteq\Fr$, $\Cn{X}$ is called the {\cmS}-\textit{theory generated by} $X$. A subset $X\subseteq\Fr$, as well as the theory $\Cn{X}$,  is called \textit{inconsistent} if $\Cn{X}=\Fr$; otherwise both are \textit{consistent}. Thus, given a system $\deduc$, $\Th$ is one of the system's theories; that is to say, if $X\subseteq\Th$ and $X\deduc A$, then $A\in\Th$. This simple observation sheds light on the central idea of Lindenbaum method, which will be explained soon. For now, let us fix the ordered pair $\left<\Alg,\Th\right>$ and call it a Lindenbaum matrix. (The full definition will be given later.) We note that an operator \textit{Cn} satisfying $(c_1)-(c_3)$ can be obtained from a \textit{closure system} over $\Fr$; that is for any subset $\mathcal{A}\subseteq \mathcal{P}(\Fr)$, which is closed under arbitrary intersection, we define:
 \[
 \textit{Cn}_{\mathcal{A}}(X)=\cap\set{Y}{X\subseteq Y\mbox{ and }Y\in\mathcal{A}}.
 \]
 It is well known that any consequence operator can be defined in this way. (Cf.~\cite{woj88}, section 1.2.)\\

\noindent Another way of defining deductive systems is through the use of logical matrices. Given a language {\Lan}, a \textit{logical} {\Lan}-\textit{matrix} (or simply a \textit{matrix}) is a pair $\Mat=\left<\fA,\Fil\right>$, where {\fA} is an {\Lan}-algebra and $\Fil\subseteq|\fA|$, where the latter is the universe of {\fA}. The (nonempty) set {\Fil} is called a \textit{filter} of the matrix $\Mat$ and the elements of $\Fil$ are called \textit{designated}. Given a matrix $\Mat=\left<\fA,\Fil\right>$, the cardinality of $|\fA|$ is also the \textit{cardinality} of $\Mat$.\\

\noindent Given a matrix $\Mat=\left<\fA,\Fil\right>$, any homomorphism of $\Alg$ into {\fA} is called a \textit{valuation} (or an \textit{assignment}). Each such homomorphism can be obtained simply by assigning elements of $|\fA|$ to the variables of $\VarP$, since, by virtue of Theorem~\ref{P:absolutely-free}, any $v:\VarP\longrightarrow|\fA|$ can be extended uniquely to a homomorphism $\hat{v}: \Alg\longrightarrow \fA$. Usually, $v$ is meant under a valuation (or an assignment) of variables in a matrix.\\

\noindent Now let $\sigma$ be a substitution and $v$ be any assignment in an algebra {\fA}. Then, defining
\begin{equation}\label{E:substitution}
v_{\sigma}=v\circ\sigma,
\end{equation}
we observe that $v_{\sigma}$ is also an assignment in {\fA}.\\

\noindent With each matrix $\Mat=\left<\fA,\Fil\right>$, we associate a relation $\modelM$ between subsets of {\Fr} and formulas of {\Fr}. Namely we define
\[
\begin{array}{rl}
X\modelM A \Longleftrightarrow
&\mbox{for every assignment $v$, if $v[X]\subseteq\Fil$,}\\
&\mbox{then $v(A)\in\Fil$}.
\end{array}
\]
Then, we observe that the following properties hold:
\[
\begin{array}{cl}
(m_1) &\mbox{if $A\in X$, then $X\modelM A$};\\
(m_2) &\mbox{if $X\modelM B$ for all $B\in Y$, and $Y\modelM A$, then $X\modelM A$}.
\end{array}
\]
Also, with help of the definition (\ref{E:substitution}), we derive the following:
\[
\begin{array}{cl}
(m_3) &\mbox{if $X\modelM A$, then for every substitution $\sigma$, $\sigma[X]\modelM\sigma(A)$}.
\end{array}
\]

\noindent Comparing the condition $(m_1)-(m_3)$ with $(s_1)-(s_3)$, we conclude that every matrix defines a structural deductive system and hence, in view of (\ref{E:cn-system}), a structural consequence operator.\\

\noindent Given a system {\cmS}, suppose a matrix $\Mat=\left<\fA,\Fil\right>$ satisfies the condition
\begin{equation}\label{E:s-filter}
\mbox{if $X\deduc A$ and $v[X]\subseteq\Fil$, then $v(A)\in\Fil$}.
\end{equation}
Then the filter $\mathcal F$ is called an {\cmS}-\textit{filter} and the matrix {\Mat} is called an
{\cmS}-\textit{matrix} (or an {\cmS}-\textit{model}). In view of (\ref{E:s-filter}), {\cmS}-matrices are an important tool in showing that $X\deduc A$ does not hold. This idea has been employed in proving that one axiom is independent from a group of others in the search for an independent axiomatic system, as well as for semantic completeness results.\\

\noindent As Lindenbaum's famous theorem below explains, every structural system {\cmS} has an {\cmS}-model.

\begin{theo}[Lindenbaum]\label{P:lindenbaum}
For any structural deductive system $\cmS$, the matrix
$\left<\Fr,\Cn{\emptyset}\right>$ is an $\cmS$-model. Moreover, for any formula $A$,
\[
A\in\Th\Longleftrightarrow\mbox{ $v(A)\in\Cn{\emptyset}$ for any valuation $v$}.
\]
\end{theo}

\noindent A matrix $\left<\fA,\Fil\right>$ is said to be \textit{weakly adequate} for a deductive system {\cmS} if
for any formula $A$,
\[
A\in\Th\Longleftrightarrow\mbox{ $v(A)\in\Fil$ for any valuation $v$}.
\]

\noindent Thus, according to Theorem~\ref{P:lindenbaum}, every structural system {\cmS} has a weakly adequate {\cmS}-matrix of cardinality less than or equal to $\overline{\overline{\Var}}+\aleph_{0}$. In general, in the last assessment, $\aleph_0$ cannot be omitted. For instance, if $\cmS=\textit{IPC}$ (intuitionistic propositional calculus), $\cmS$ has no finite weakly adequate matrix. (Cf.~\cite{god32}.)\\

\noindent An {\cmS}-matrix is called \textit{ strongly adequate} for {\cmS} if for any set $X\subseteq\Fr$ and any formula $A$,
\begin{equation}\label{E:adequate}
X\deduc A\Longleftrightarrow X\modelM A.
\end{equation}

\noindent We note that Theorem~\ref{P:lindenbaum} cannot be improved to include strong adequacy. Also,
if $\overline{\overline{\Var}}\leq\aleph_{0}$ and $\cmS=\textit{IPC}$, there is no denumerable matrix {\Mat} with (\ref{E:adequate}). (Cf.~\cite{wro74a}.)\\

\noindent\textbf{Historical remarks}\\ 
\noindent A. Tarski seems to be the first who promoted ``the view of matrix formation as a general method of constructing systems''~\cite{lt30}. However, matrices had been employed earlier, e.g., by P. Bernays~\cite{ber26} and others either in the search for an independent axiomatic system or for defining a system different from classical logic. Also, later on J.C.C. McKinsey~\cite{mck39} used matrices to prove independence of logical connectives in intuitionistic propositional logic.\\

\noindent Theorem~\ref{P:lindenbaum} was discovered by A. Lindenbaum. Although this theorem was not published by the author, it had been known in Warsaw-Lvov logic circles at the time. In a published form it appeared for the first time in~\cite{lt30} without proof. Its proof appeared later on in the two independent publications~\cite{los49} and~\cite{her51}. McKinsey and Tarski~\cite{mt48} gave an example of a deductive system with $\overline{\overline{\Var}}\leq\aleph_{0}$ but without any finite weakly adequate matrix.\\

\noindent\Large{\textbf{W\'{o}jcicki's Theorems}}\\

\noindent We get more {\cmS}-matrices, noticing the following. Let {\Thy} be an {\cmS}-theory. The pair  $\left<\Fr,\Thy\right>$ is called a \textit{Lindenbaum matrix} relative to {\cmS}. We observe that  for any substitution $\sigma$,
\begin{equation*}\label{E:-theory-as-filter}
\mbox{if  $X\deduc A$  and $\sigma[X]\subseteq\Thy$, then $\sigma(A)\in\Thy$}.
\end{equation*}
That is to say,  any Lindenbaum matrix relative to a system {\cmS} is an {\cmS}-model.\\

\noindent A deductive system {\cmS} is said to be \textit{uniform} if, given a set $X\subseteq\Fr$ and a consistent set $Y\subseteq\Fr$, $X\cup Y\deduc A$ and $\Var(Y)\cap \Var(A)=\emptyset$ imply $X\deduc A$. A system {\cmS} is \textit{couniform} if for any collection $\{X_{i}\}_{i\in I}$ of formulas with $\Var(X_{i})\cap \Var(X_{j})=\emptyset$, providing $i\neq j$, if the set $\cup\{X_{i}\}_{i\in I}$ is inconsistent, then at least one $X_{i}$ is inconsistent as well.

\begin{theo}[W\'{o}jcicki]\label{P:wojcicli-1}
A structural deductive system $\cmS$ has a strongly adequate matrix if and only if {\cmS} is both uniform and couniform.
\end{theo}

\noindent For the ``if'' implication of the statement, the matrix of Theorem~\ref{P:lindenbaum} is not enough. However, it is possible to extend the original language {\Lan} to $\Lan^{+}$ in such a way that the \textit{natural extension} $\textit{Cn}_{\mathcal{S}^{+}}$ of $\textit{Cn}_{\mathcal{S}}$ onto $\Lan^{+}$ allows one to define a Lindenbaum matrix $\left<\mathfrak{F}_{\mathcal{L}^{+}},\textit{Cn}_{\mathcal{S}^{+}}(X)\right>$, for some $X\subseteq\textit{Fr}_{\mathcal{L}^{+}}$, which is  strongly adequate for {\cmS}. (Cf.~\cite{woj88} for detail.)\\

\noindent A pair $\left<\fA,\{\Fil_{i}\}_{i\in I}\right>$,  where  {\fA} is an {\Lan}-algebra and each $\Fil_{i}\subseteq|\fA|$, is called a \textit{generalized matrix} (or a {\g}-\textit{matrix} for short). A \g-matrix is a \g-\cmS-\textit{model} (or a \g-\cmS-\textit{matrix}) if each $\left<\fA,\Fil_{i}\right>$ is an \cmS-model. (In~\cite{dh01} a \g-matrix is called an \textit{atlas}.)

\begin{theo}[W\'{o}jcicki]
For every structural deductive system {\cmS}, there is a \g-\cmS-matrix {\Mat} of cardinality $\overline{\overline{\Var}}+\aleph_{0}$, which is strongly adequate for {\cmS}.
\end{theo}

\noindent Indeed, let $\{\Thy\}$ be the collection of all \cmS-theories. Then the \g-matrix $\left<\Fr,\{\Thy\}\right>$ is strongly adequate for {\cmS}. (Cf.~\cite{woj88},~\cite{dh01} for detail.)\\

\noindent We note that, alternatively, one could use the notion of a \textit{bundle of matrices}; a bundle is a set $\set{\left<\fA,\Fil_{i}\right>}{i\in I}$, where {\fA} is an \Lan-algebra and each $\Fil_{i}$ is a filter of {\fA}. (Cf.~\cite{woj88}, section 3.2.11.)\\

\noindent\textbf{Historical remarks}\\ Theorem~\ref{P:wojcicli-1} was the result of the correction by R. W\'{o}jcicki of an erroneous assertion in~\cite{ls58}, where the important question on the strong adequacy of a system was raised. A number of algebraic equivalents of uniformity is discussed in~\cite{dkw09}.\\

\noindent T. Smiley~\cite{smi62} was perhaps the first to propose \g-matrices (known also as \textit{Smiley matrices}) defined as pairs $\left<\fA,\textit{Cn}\right>$, where {\fA} is an {\Lan}-algebra and an operator $\textit{Cn}:\mathcal{P}(|\fA|)\rightarrow\mathcal{P}(|\fA|)$ satisfies the conditions $(c_1)-(c_3)$ (with \textit{Cn} instead if $\textit{Cn}_{\mathcal{S}}$). Then, Smiley defined $x_1,\ldots, x_n\vdash y$ if and only of $y\in \textit{Cn}(\{x_1,\ldots,x_n\})$, where it is assumed that $|\fA|\subseteq U$, where $U$ is a universal set of sentences.\\

\noindent\Large{\textbf{Lindenbaum-Tarski Algebra}}\\

\noindent The question of the possibility to decide, whether $X\deduc A$ is true or not is central in theory of deduction. Although the notion we are about to introduce is less general than that of {\cmS}-matrix, it points out at a way, following which this question can be often fruitfully discussed.\\

\noindent An {\cmS}-matrix $\left<\fA,\Fil\right>$ is said to be \textit{univalent} (or an {\Sr}-\textit{matrix}) if the {\cmS}-filter $\Fil$ consists of one value, say $\Fil=\{\on\}$, where $\on\in|\fA|$. Let us restrict our original question to the following: How can the property $\emptyset\deduc A$ be characterized in matrix terms?

Let $\left<\fA,\{\on\}\right>$ be an {\Sr}-matrix and $A$ be an {\cmS}-theorem. Then, in view of (\ref{E:s-filter}), $v(A)=\on$ for every valuation $v$ in {\fA}. It would be interesting to know when the converse is true too. Thus the main problem is: How can one obtain an {\Sr}-matrix?

\begin{defi}[Lindenbaum-Tarski algebra]\label{D:lindenbaum-alg}
Let {\Thy} be an {\cmS}-theory and let $\Theta(\Thy)$ be the congruence on {\Alg} generated by {\Thy}; \emph{cf.~\cite{bs81}}. The quotient algebra $\Alg/\Theta(\Thy)$ is called a Lindenbaum-Tarski algebra of {\cmS} relative to {\Thy}. If $\Thy=\Th$, then we call this quotient simply a Lindenbaum-Tarski algebra.
\end{defi}

\noindent An important conclusion from this definition is the following.

\begin{theo}\label{P:lindenbaum-tarski}
Let {\cmS} be a structural deductive system and {\Thy} be a nonempty {\cmS}-theory. Assume that {\Thy} is a congruence class with respect to $\Theta(\Thy)$. Then $\left<\Alg/\Theta(\Thy),\{\Thy\}\right>$ is an {\Sr}-matrix; that is to say, denoting $\on=\Thy$, if $X\deduc A$ and $v$ is a valuation in $\Alg/\Theta(\Thy)$, then
\begin{equation}\label{E:lindenbaum-general}
v[X]=\{\on\} \Longrightarrow v(A)=\on.
\end{equation}
Moreover, if $\Thy=\Th$, then
\begin{equation}\label{E:lindenbaum-simple}
A\in\Th \Longleftrightarrow \mbox{$v(A)=\on$ for any valuation $v$ in $\Alg/\Theta(\Th)$}.
\end{equation}
Let the valuation $v_{0}(p)=p/\Theta(\Th)$ for every $p\in\Var$. Then
\begin{equation}\label{E:lindenbaum-simple2}
A\in\Th \Longleftrightarrow v_{0}(A)=\on.
\end{equation}
\end{theo}

\begin{defi}
Let {\cmS} be a structural deductive system. We say that  {\cmS} admits the Lindenbaum-Tarski algebra $($relative to {\Thy}$)$ if {\Th} $(${\Thy} respectively$)$ is a congruence class with respect to $\Theta(\Th)$ $($with respect to $\Theta(\Thy)$$)$ on {\Alg}.
\end{defi}

%Then the question arises: Which properties should %{\deduc} have in order $\Thy$ to be a congruence class %with respect to $\Theta(\Thy)$?

\noindent Now let us convert the propositional language $\Lan$ into
a first order language $\LanP$ with equality so that the propositional variables  and the logical connectives of $\Lan$ become the individual variables and functional constants of $\LanP$, respectively. The set of individual variables is denoted by $\VarPP$.  Also, $\LanP$ has an individual constant $\on$,  the equality symbol `$=$' and universal and existential quantifiers. (Actually, we will need only the former.) We can assume that there is no logical connectives in $\LanP$. Since the formulas of $\Lan$ now become terms of $\LanP$, each \textit{atomic formula} of $\LanP$ is an expression of the form:
\begin{center}
\begin{tabular}{l}
$A(p,\on,\ldots)=B(q,\on,\ldots)$,
\end{tabular}
\end{center}
where variables $p$ and $q$ are not necessarily distinct and they, as well as the constant $\on$, may or may not occur in the equality.\\

\noindent A universal closure (in the sense of first order logic) of an atomic formula of $\LanP$ is often referred to as an \textit{identity}. We will deal with interpretations of identities only. Therefore, we semantically treat atomic formulas and their universal closures equally. An unspecified identity will be denoted by $\varphi$.\\

\noindent The $\LanP$-formulas are interpreted in algebras {\fB} of the type $\Lan$ endowed with a 0-ary operation $\on$. Then, for instance, an identity
\[
A(p,\on,\ldots)=\on
\]
is said to be \textit{valid} (or to \textit{hold}) in $\fB$, in symbols $\fB\models A(p,\on,\ldots)=\on$, if for any assignment $v:\VarP\rightarrow |\fB|$
\[
A(v(p),\on,\ldots)=\on.
\]

\noindent Given a system {\cmS}, we denote
\[
\AlgS=\left<\Alg/\Theta(\Th), \on\right>,
\]
where $\on$ is the congruence class generated by {\Th}. Thus $\AlgS$ is the expansion of $\Alg/\Theta(\Th)$ obtained by adding the constant $\on$ to the signature of the latter.
Then, we define:
\[
\Phi_{\mathcal{S}}=\set{A=\on}{A\in\Th}
\]
and
\[
K_{\mathcal{S}}=\set{\fB}{\fB\models\varphi\mbox{ for all }\varphi\in\Phi_{\mathcal{S}}}.
\]
It is obvious that the class $K_{\mathcal{S}}$ is a variety.

\begin{theo}\label{P:lindenbaum-allowed}
Let a structural deductive system {\cmS} admit the Lindenbaum-Tarski algebra. Then the algebra $\AlgS$ belongs in the variety $K_{\mathcal{S}}$. More than that,
\[
\AlgS\models A=\on\Longleftrightarrow A\in\Th.
\]
Moreover, $\AlgS\models A(p,\on,\ldots)=\on$ if and only if $A(p/\Theta(\Th),\on,\ldots)=\on$ in $\AlgS$, that is
$A(p/\Theta(\Th),\Th,\ldots)=\Th$ in $\Alg/\Theta(\Th)$.
\end{theo}

\noindent Theorem~\ref{P:lindenbaum-allowed} gives rise to the following questions: When is $\AlgS$ functionally free~\cite{tar46} in $K_{\mathcal{S}}$? When is $\AlgS$ a free algebra in $K_{\mathcal{S}}$?\\

\noindent\textbf{Historical remarks}\\ 
\noindent In two  parts,~\cite{tar35} and~\cite{tar36}, of one paper, the English translation of which constitutes one chapter, \textit{Foundations of the Calculus of Systems}, of~\cite{tar56}, A. Tarski showed that the Lindenbaum-Taski algebra of the system based on classical propositional calculus is a Boolean algebra.\\

\noindent\Large{\textbf{Alternative Approach}}\\

\noindent Let $\left<\fA,\mathcal{F}\right>$ be a matrix. A congruence (or an equivalence) $\theta$ on {\fA} is said to be \textit{compatible} with $\mathcal{F}$ if $\cup\{ x /\theta ~|~x\in\mathcal{F}\}=\mathcal{F}$. Since the identity relation is compatible with any $\mathcal{F}$, the set of compatible congruences (or equivalences) is not empty for any matrix. Then, it can be proven~\cite{bp89} that for any matrix $\Mat=\left<\fA,\mathcal{F}\right>$, there is a largest congruence of {\fA} compatible with $\mathcal{F}$. This congruence is called the \textit{Leibniz congruence} of {\Mat}; it is denoted by $\Omega_{\mathfrak{A}}\mathcal{F}$ and can be defined as follows:
\[
\begin{array}{rl}
\Omega_{\mathfrak{A}}{\mathcal{F}}=
&\lbrace(a,b)~|~\forall
A(p,p_{0},\ldots,p_{n})\forall c_{0},\ldots, c_{n}\in|\mathfrak{A}|.\\
&A(a,c_{0},\ldots,c_{n})\in{\mathcal{F}}\Leftrightarrow A(b,c_{0},\ldots,c_{n})\in{\mathcal{F}}\rbrace.
\end{array}
\]

\noindent If the matrix in question is a Lindenbaum one, say $\left<\Alg,\Thy\right>$, then an example of a compatible equivalence on this matrix is a \textit{Frege relation} $ \Lambda\Thy $ defined as follows:
\[
(A,B)\in\Lambda\Thy \Longleftrightarrow \Thy, A\deduc B\mbox{ and }\Thy, B\deduc A\tag{\textit{Frege relation relative to} $\Thy$}
\]

\vspace{0.15in}
\noindent A system {\cmS} is called \textit{Fregean} if each $\Lambda\Thy$ is a congruence on {\Alg}. Obviously, if {\cmS} is Fregean, it admits the Lindenbaum-Tarski algebra relative to any $\Thy$.\\

\noindent Another example of a compatible relation on $\left<\Alg,\Thy\right>$ is the largest congruence of {\Alg} contained in $\Lambda\Thy$, which is referred to as a \textit{Suszko congruence}:
\[
\begin{array}{rl}
(A,B)\in\Tilde{\Omega}\Thy\Longleftrightarrow
&\mbox{for every $C(p)$, $\Thy, C(A/p)\deduc C(B/p)$}\\
&\mbox{and $\Thy,C(B/p)\deduc C(A/p)$}.
\end{array}
 \tag{\textit{Suszko congruence relative to} $\Thy$}
\]

\vspace{0.15in}
\noindent Obviously, a system {\cmS} is Fregean if and only if $\Lambda\Thy=\Tilde{\Omega}\Thy$ for all {\Thy}.\\

\noindent The Leibniz congruence of a matrix $\left<\Alg,\Thy\right>$ is referred to as \textit{Leinbniz congruence relative to} {\Thy}. It  turns out that
\[
\Omega\Thy=\cap\set{\Tilde{\Omega}\Sigma^{\prime}_{\mathcal{S}}}{\Thy\subseteq\Sigma^{\prime}_{\mathcal{S}}}
\]
and, therefore, each Suszko congruence $\Tilde{\Omega}\Thy$ is compatible with $\Thy$. Also, given a system $\mathcal{S}$, one defines
\[
\Tilde{\Omega}_{\mathcal{S}}=\cap\{\Tilde{\Omega}\Sigma_{\mathcal{S}}~|~
\Sigma_{\mathcal{S}}\mbox{ is an $\mathcal{S}$-theory}\}.\tag{\emph{Tarski congruence}}
\]

\noindent Thus we have:
\[
\Tilde{\Omega}_{\mathcal{S}}\subseteq\Tilde{\Omega}\Thy\subseteq\Lambda
\Thy\cap\Omega\Thy.
\]

\noindent Suszko, Leibniz and Tarski congruences give rise to the {\cmS}-matrices $\left<\Alg/\Omega\Thy,\Thy/\Omega\Thy\right>$, $\left<\Alg/\Tilde{\Omega}\Thy,\Thy/\Tilde{\Omega}\Thy\right>$, and the \emph{g}-$\mathcal{S}$-matrix $\left<\Alg/\Tilde{\Omega}_{\mathcal{S}},\{\Thy/\Tilde{\Omega}_{\mathcal{S}}
~|~\Thy\mbox{ is an $\mathcal{S}$-theory}\}\right>$, whose first components, $ \Alg/\Omega\Thy $,
$ \Alg/\Tilde{\Omega}\Thy $ and $ \Alg/\Tilde{\Omega}_{\mathcal{S}} $, in \textit{Algebraic abstract logic} are also referred to as \textit{Lindenbaum-Tarski algebras}. (See~\cite{fjp03} and~\cite{fj09} for comprehensive surveys.)\\

\noindent\Large{\textbf{Specifications and Applications}}\\

\noindent A structural deductive system  {\cmS} is called \textit{implicative extensional} if
its  language {\Lan} contains a binary connective $ \rightarrow $ (will be written in the infix notation), and for any {\cmS}-theory $\Thy$ and any $A, B, C\in\Fr$, the following conditions hold:
\[
\begin{array}{cl}
(i_1) & A\rightarrow A\in\Thy;\\
(i_2) &B\in\Thy\Longrightarrow A\rightarrow B\in\Thy;\\
(i_3) &A\rightarrow B, B\rightarrow C\in\Thy \Longrightarrow
A\rightarrow C\in\Thy;\\
(i_4) &A, A\rightarrow B\in\Thy \Longrightarrow B\in\Thy;\\
(i_5) &A_{i}\rightarrow B_{i}, B_{i}\rightarrow A_{i}\in\Thy, 1\leq i\leq n, \Longrightarrow \Pi A_{1}\ldots A_{n}\rightarrow\Pi B_{1}\ldots B_{n}\\
&\mbox{for each $n$-ary connective $\Pi$}.
\end{array}
\]

\noindent Now, given {\cmS}, we consider the following relation on $\Fr$:
\[
A\approx_{\mathcal{S}} B\Longleftrightarrow A\rightarrow B, B\rightarrow A\in\Th.\tag{\textit{Rasiowa relation}}
\]
\begin{theo}[Rasiowa]\label{P:rasiowa}
If {\cmS} is an implicative extensional system, then the relation $\approx_{\mathcal{S}}$ is a congruence on {\Alg}. Moreover, $\Th$ is a congruence class with respect to $\approx_{\mathcal{S}}$.
\end{theo}

\noindent Applying Theorem~\ref{P:rasiowa}  to \textit{IPC}, one can observe (actually, it was shown in~\cite{rs70}) that
$\Alg/\!\!\approx_{\textit{IPC}}$ is the free algebra of rank $\overline{\overline{\Var}}$ in the variety of Heyting algebras. Using the Tarski relation $\approx_{\textit{IPC}} $, Nishimura~\cite{nis60} gave an elegant description of the Lindenbaum-Tarski algebra of \textit{IPC} in a language with a single propositional variable. This algebra is also the free algebra of rank $1$ in the variety of Heyting algebras. See \textit{Free algebra}.\\

\noindent Also, it is worth noticing that, using a Lindenbaum-Tarski algebra as defined above, one can prove that there is an algorithm which decides whether two finite \g-matrices define the same deductive system; this result is due to A. Citkin (unpublished) and J. Zygmunt~\cite{zyg83}. In this connection see \textit{Decision problem}.\\

\noindent\textbf{Historical remarks}\\
\noindent In~\cite{tar35},~\cite{tar36} (see~\cite{tar56}, chapter XII), Tarski gave the first specification of a system which admits the Lindenbaum-Tarski algebra. Later on, Rasiowa~\cite{ras74} summarized the work that had been done by the time in the notion of ``the class of standard systems of implicative extensional propositional calculi,'' which is a simplified version of that we use above.\\

\noindent Also, if {\cmS} is an implicative extensional system, then $\AlgS$ as defined above is Rasiowa's {\cS}-algebra~\cite{ras74}, or nowadays known~\cite{bp89} as \textit{Hilbert algebra} with compatible operations.\\

\noindent In~\cite{zyg83} Zygmunt credits Citkin for the decidability result mentioned above. Recently, it was rediscovered by L. Devyatkin~\cite{dev13}.

%\bibliographystyle{amsplain}
%\bibliography{Bibtex}

\begin{thebibliography}{99.}
	
\bibitem{ber26}
Paul Bernays, \emph{Untersuchung des {A}ussagenkalk\"{u}ls der ``{P}rincipia
	{M}athematica''}, Math. Z. \textbf{25} (1926), 305--320. (German)

\bibitem{bp89}
W.~J. Blok and Don Pigozzi, \emph{Algebraizable Logics}, Mem. Amer. Math. Soc.
\textbf{77} (1989), no.~396, vi+78. \MR{973361 (90d:03140)}

\bibitem{bs81}
Stanley Burris and H.~P. Sankappanavar, \emph{A Course in Universal Algebra},
Graduate Texts in Mathematics, vol.~78, Springer-Verlag, New York, 1981.
\MR{648287 (83k:08001)}

\bibitem{dev13}
Leonid~Yu. Devyatkin, \emph{Equality of consequence relations in finite-valued logical matrices}, Logical Investigations, vol.~19, 273--280, Ross. Akad. Nauk, Inst. Filos., Moscow, 2013. \MR{3204258}

\bibitem{dh01}
J.~Michael Dunn and Gary~M. Hardegree, \emph{Algebraic Methods in Philosophical
	Logic}, Oxford Logic Guides, vol.~41, The Clarendon Press Oxford University
Press, New York, 2001, Oxford Science Publications. \MR{1858927
	(2002j:03001)}

\bibitem{dkw09}
W.~Dziobiak, A.~V. Kravchenko and P.~J. Wojciechowski,
\emph{Equivalents for a quasivariety to be generated by a single structure}, Studia Logica \textbf{91} (2009), no.~1, 113--123. \MR{2476339 (2010a:08010)}

\bibitem{fjp03}
J.~M. Font, R.~Jansana, and D.~Pigozzi, \emph{A survey of abstract algebraic
	logic}, Studia Logica \textbf{74} (2003), no.~1-2, 13--97, Abstract algebraic
logic, Part II (Barcelona, 1997). \MR{1996593 (2004m:03241)}

\bibitem{fj09}
Josep~Maria Font and Ramon Jansana, \emph{A General Algebraic Semantics for
	Sentential Logics}, second edition, Lecture Notes in Logic, vol.~7.

\bibitem{god32}
Kurt G\"{o}del, \emph{Zum intuitionistischen {A}ussagenkalk\"{u}l}, Anzeiger der {A}kademie der {W}issenschaften, mathematisch-naturwissenschaftiche {Kl}lasse \textbf{69} (1932), 65--66; (German) reprinted in:~\cite{god86}, p. 224. English edition in:~\cite{god86}, pp. 223 and 225.

\bibitem{god86}
Kurt G\"{o}del, \emph{Collected Works, Vol. 1, Publications 1929--1936}, S. Feferman (editor), The Clarendon Press, Oxford University Press, New York, 1986.

\bibitem{her51}
Hans Hermes, \emph{Zur {T}heorie der aussagenlogischen {M}atrizen}, Math. Z.
\textbf{53} (1951), 414--418. (German) \MR{0040241 (12,663c)}

\bibitem{ls58}
J.~{\L}o{\'s} and R.~Suszko, \emph{Remarks on sentential logics}, Nederl. Akad.
Wetensch. Proc. Ser. A 61 = Indag. Math. \textbf{20} (1958), 177--183.
\MR{0098670 (20 \#5125)}

\bibitem{los49}
Jerzy {\L}o{\'s}, \emph{On logical matrices}, Trav. Soc. Sci. Lett. Wroc\l aw.
Ser. B. \textbf{1949} (1949), no.~19, 42. \MR{0089812 (19,724b)}

\bibitem{lt30}
Jan {\L}ukasiewicz and Alfred Tarski, \emph{Untersuchungen \"{u}ber den
	{A}ussagenkalk\"{u}l}, Comptes rendus des s\'{e}ances de la Soci\'{e}t\'{e}
des Sciences et des Lettres de Varsovie, \textbf{23} (1930), cl. iii, 30--50. (German)
English edition in:~\cite{tar56}, pp. 38--59.

\bibitem{mck39}
J.~C.~C. McKinsey, \emph{Proof of the independence of the primitive symbols of
	{H}eyting's calculus of propositions}, J. Symbolic Logic \textbf{4} (1939),
155--158. \MR{0000805 (1,131f)}

\bibitem{mt48}
J.~C.~C. McKinsey and Alfred Tarski, \emph{Some theorems about the sentential calculi of {L}ewis and
{H}eyting}, J. Symbolic Logic \textbf{13} (1948),
1--15. \MR{0024396 (9,486q)}

\bibitem{mmt87}
Ralph~N. McKenzie, George~F. McNulty and Walter~F. Taylor, \emph{Algebras, Lattices, Varieties}, vol.~1,
Wadsworth \& Brooks/Cole Advanced Books \& Software, Monterey, CA, 1987. \MR{883644 (88e:08001)}

\bibitem{nis60}
Iwao Nishimura, \emph{On formulas of one variable in intuitionistic
	propositional calculus.}, J. Symbolic Logic \textbf{25} (1960), 327--331
(1962). \MR{0142456 (26 \#25)}

\bibitem{ras74}
Helena Rasiowa, \emph{An Algebraic Approach to Non-Classical Logics},
North-Holland Publishing Co., Amsterdam, 1974, Studies in Logic and the
Foundations of Mathematics, Vol. 78. \MR{0446968 (56 \#5285)}

\bibitem{rs70}
Helena Rasiowa and Roman Sikorski, \emph{The Mathematics of Metamathematics},
third ed., PWN---Polish Scientific Publishers, Warsaw, 1970, Monografie
Matematyczne, Tom 41. \MR{0344067 (49 \#8807)}

\bibitem{smi62}
Timothy Smiley, \emph{The independence of connectives}, J. Symbolic Logic
\textbf{27} (1962), 426--436. \MR{0172784 (30 \#3003)}

\bibitem{sur82}
Stanis{\l}aw~J. Surma, \emph{On the origin and subsequent applications of the
	concept of the {L}indenbaum algebra}, Logic, methodology and philosophy of
science, {VI} ({H}annover, 1979), Stud. Logic Foundations Math., vol. 104,
North-Holland, Amsterdam, 1982, pp.~719--734. \MR{682440 (84g:01045)}

\bibitem{tar35}
Alfred Tarski, \emph{Grundz\"{u}ge der Systemenkalk\"{u}l. Erster Teil},
Fundamenta Mathematica \textbf{25} (1935), 503--526.

\bibitem{tar36}
Alfred Tarski, \emph{Grundz\"{u}ge der Systemenkalk\"{u}l. Zweiter Teil}, Fundamenta
Mathematica \textbf{26} (1936), 283--301. (German)

\bibitem{tar46}
Alfred Tarski, \emph{A remark on functionally free algebras}, Ann. of Math. (2)
\textbf{47} (1946), 163--165. \MR{0015038 (7,360a)}

\bibitem{tar56}
Alfred Tarski, \emph{Logic, {S}emantics, {M}etamathematics. {P}apers from 1923 to
	1938}, Oxford at the Clarendon Press, 1956, Translated by J. H. Woodger.
\MR{0078296 (17,1171a)}

\bibitem{woj88}
Ryszard W{\'o}jcicki, \emph{Theory of Logical Calculi}, Synthese Library, vol.
199, Kluwer Academic Publishers Group, Dordrecht, 1988, Basic theory of
consequence operations. \MR{1009788 (90j:03001)}

\bibitem{wro74a}
Andrzej Wro{\'n}ski, \emph{On cardinalities of matrices strongly adequate for
	the intuitionistic propositional logic}, Rep. Math. Logic (1974), no.~3,
67--72. \MR{0387011 (52 \#7858)}

\bibitem{zyg83}
Jan Zygmunt, \emph{An application of {L}indenbaum method in the domain of strongly finite sentential calculi}, Prace Foilozoficzne, \textbf{29} (1983), 59--68.

\end{thebibliography}

\stop
\subsection*{Algebarizable deductive systems}

Some definitions.\\
Given a deductive system {\cmS},
\[
(A,B)\in\Lambda{\mathcal{S}}\Longleftrightarrow A\deduc B\mbox{ and }B\deduc A.\tag{\textit{Frege relation of} {\cmS}}
\]

A system {\cmS} is called \textit{selfextensional}~\cite{woj88}, Chapter 5, if $\Lambda\cmS$ is a congruence on {\Alg}.

Given an {\cmS}-theory {\Thy},
\[
(A,B)\in\Lambda\Thy \Longleftrightarrow \Thy, A\deduc B\mbox{ and }\Thy, B\deduc A\tag{\textit{Frege relation relative to} $\Thy$}
\]
\begin{prop}
\[
(A,B)\in\Lambda\Thy \Longleftrightarrow\forall\mathcal{S}\mbox{-theory }\Sigma.
~\Thy\subseteq\Sigma\Rightarrow [A\in\Sigma\Leftrightarrow
B\in\Sigma].
\]
\end{prop}

It is obvious that $\Lambda\Thy$ is reflexive and transitive.
A system {\cmS} is called \textit{Fregean} (or \textit{extensional}) if for any theory {\Thy}, $\Lambda\Thy$ is a congruence on {\Alg}.
\begin{cor}
If $\Lambda\Thy$ is a congruence on {\Alg}, then
\[
A\in\Thy\Longrightarrow A/\Lambda\Thy\subseteq\Thy.
\]
\end{cor}

Regardless whether $\Lambda\Thy$ is a congruence or not, there is a largest congruence on {\Alg}, which is contained in $\Lambda\Thy$. This congruence was introduced by Suszko  and is denoted by $\Tilde{\Omega}\Thy$. It can be defined as follows: Given {\cmS}-theory $\Thy$,
\[
(A,B)\in\Tilde{\Omega}\Thy\Longleftrightarrow\mbox{ for every $C(p)$, $C(A/p)\deduc C(B/p)$ and $\Thy,C(B/p)\deduc C(A/p)$}.\tag{\textit{Suszko congruence relative to} $\Thy$}
\]

\begin{cor}
\[
A\in\Thy\Longrightarrow A/\Tilde{\Omega}\Thy\subseteq\Thy.
\]
\end{cor}

It can be proven that
\[
(A,B)\in\Tilde{\Omega}\Thy\Longleftrightarrow\mbox{ for every $C(p)$ and $\Sigma^{\prime}_{\mathcal{S}}\supseteq\Thy$, $C(A/p)\in
\Sigma^{\prime}_{\mathcal{S}}\Leftrightarrow C(B/p)\in\Thy^{\prime}$}.
\]

Obviously, {\cmS} is Fregean if and only if $\Tilde{\Omega}\Thy=\Lambda\Thy$ for all {\Thy}.

On the other hand, we have these definitions.

Let $\left<\fA,\mathcal{F}\right>$ be a matrix. A congruence $\theta$ is \textit{compatible} with $\mathcal{F}$ if $\cup\{ x /\theta ~|~x\in\mathcal{F}\}=\mathcal{F}$. Since the identity relation is compatible with any $\mathcal{F}$, the set of compatible congruences is not empty for any matrix. Then, it can be proven that for any matrix $\Mat=\left<\fA,\mathcal{F}\right>$, there is a largest congruence on {\fA} compatible with $\mathcal{F}$. We call it a \textit{Leibniz congruence} and denote it by $\Omega_{\mathfrak{A}}\mathcal{F}$.

\section*{Preparation}
\begin{enumerate}
\item Block 1:
\begin{itemize}
\item ``The view of matrix formation as a general method of constructing systems is due to Tarski'' (cf.~\cite{lt30}; English translation in~\cite{tar56}, chapter IV; p. 40)
\item Definition of a logical (normal) matrix (cf.~\cite{lt30}, Definition 3; see~\cite{tar56}, p. 41) is due to Tarski.
\item The term \emph{designated elements} is due to Bernays (cf.~\cite{ber26}, p. 316).
\item Definition 1 in [\cite{tar56}, p. 33] A set $X$ is called a deductive system (or simply a system), in symbols $X\in\fC$, if $Cn(X)=X\subseteq S$.
\item Definition 1 in~\cite{lt30} (see~\cite{tar56}, p. 39): The set $S$ of all sentences is the intersection of all those sets which contain all sentential variables (elementary sentences) and are closed under the operations of forming implications and negations.
\item Definition 2 in~\cite{lt30} (see~\cite{tar56}, p. 40): The set of consequences $\cn(X)$ of the set $X$ of sentences is the intersection of those sets which include $X\subseteq S$ and are closed under the operations of substitution and detachment.
\item Definition of a (logical) matrix (due to Tarski): $\fM=(A,D, \rightarrow,\neg)$, where both $A$ and $D$ are nonempty and $D\subset A$.

    The matrix {\fM} is called normal if the following holds:
    \[
    x\in D\mbox{ and } y\in A\setminus D \Rightarrow x\rightarrow y\in A\setminus D.
    \]
     %(Definition 3 in~\cite{lt30}; the elements if $D$ %are called, according to~\cite{ber26}, p. 316, %\emph{designated}.
\item Theorem 2 in~\cite{lt30} (see~\cite{tar56}, p. 41): If $\fM$ is a normal matrix, then $\fG(\fM)\in\fC$.
\item Theorem 3 (Lindenbaum; see~\cite{lt30} or~\cite{tar56}, p. 38): For every system $X\in\fC$ there exists a normal matrix $\fM=(A, D, \rightarrow, \neg)$, with an at most denumerable set $A$, which satisfies the formula $X=\fG(\fM)$.
\item Definition 3 in~\cite{tar56}, p. 34: A set $X$ of sentences is called consistent, in symbols $X\in\fW$, if $X\subseteq S$ and if the formula $X\sim S$ does not hold (i.e. if $\cn(X)\neq S$).
\item Theorem 4 in~\cite{lt30} (see~\cite{tar56}, p. 42): If {\fM} is a normal matrix, then $\fG(\fM)\in\fW$.
\item Lemma 1.4 (\cite{pw08}, p. 4): Algebra {\cS} of a propositional language is free over the class of all similar algebras and the propositional variables are the free generators of {\cS}.
\end{itemize}

\item Block 2 (logical language):
\begin{itemize}
\item The important example of a propositional language, along with the symbols of {\Var}, has the connectives: {\con} (conjunction), {\dis} (disjunction), {\ra} (implication), and $\neg$ (negation). The arities are as follows: $\#(\con)=\#(\dis)=\#(\ra)=2$ (those are binary connectives) and $\#(\neg)=1$. The language with these connectives is called \emph{assertoric}.
\item In any assertoric language, the following agreement about formula formation is in use. For any binary connective $\medcirc$, instead of $\medcirc AB$ (prefix notation) one writes $(A\medcirc B)$ (infix notation). Another additional agreement allows one to omit the most outer parentheses. These are examples of assertoric statements: $A\ra (B\ra A)$, $(A\con B)\ra A$, and $A\ra(A\dis B)$, where the most outer parentheses are omitted. We note that parentheses `(' and `)' do not belong to the language under consideration but are symbols of metalanguage; they can be eliminated if we use the regular formula formation.
\end{itemize}
\section*{Deleted material}

One of the important questions is: How can a deductive system be defined? To address this question, we begin with the definition of a (\textit{structural}) \textit{rule of inference}. (We advise the reader that the structurality of a rule appears, when we apply the rule, while constructing a derivation.)

A (\textit{finitary}) \textit{rule of inference} is an ordered pair $\left<{\Ga},A\right>$, usually written as $\frac{\Ga}{A}$, where $\Ga\fsub\Fr$ and $A\in\Fr$. The formulas in $\Ga$ are called the \textit{premises} of the rule and $A$ is called the \textit{conclusion}.

Given a nonempty set $\cmS$ of rules, a (\textit{linear}) \textit{derivation with respect to} $\cmS$ (or simply $\cmS$-\textit{derivation}) of a formula $A$ from a subset $X\subseteq\Fr$ is a finite list of formulas with $A$ being the last formula on the list and each formula $B$ of the list is either a formula of $X$ or is obtained by an application of a rule $\frac{\Ga}{C}\in\cmS$, when, for some substitution $\sigma$, $B=\sigma(C)$ and the formulas $\sigma[\Ga]$ appear before $B$ on the list. Denoting the existence of an $\cmS$-derivation of $A$ from $X$ by $X\vdash_{\cmS}A$, we claim that the so-defined relation is a finitary deductive system, that is, it satisfies the conditions $(s_{1})-(s_{4})$ above.
The last assertion should not be surprising. For, given a finitary deductive system $\vdash_{\cmS}$, if we define the set $\cmS^{\ast}$ of rules as follows:
\[
\frac{\Ga}{A}\in\cmS^{\ast}\Longleftrightarrow \Ga\vdash_{\cmS}A,
\]
for all $\Ga\fsub\Fr$, then it is not difficult to show that the relations $\deduc$ and $\vdash_{\mathcal{S}^{\ast}}$ are equal.\\

Sometimes it will be convenient to denote a deductive system {\deduc} simply by {\cmS}.\\

 In the sequel we will not distinguish an assignment $v$ and its extension $\hat{v}$ uniquely obtained from the former.
 
 \section*{Deleted references}
 
 \bibitem{gra08}
 George Gr\"{a}tzer, \emph{Universal algebra}, Springer, second edition with updates, 2008.
 \MR{2455216 (2009j:08001)}

\end{enumerate}

\end{document}